\documentclass[12pt]{article}          \usepackage{amsfonts} \usepackage{amssymb} \usepackage{latexsym} \usepackage{graphicx} \usepackage{amsbsy} 

\newtheorem{thm}{Theorem}

\newtheorem{lemma}[thm]{Lemma}

\newcommand{\bprf}[1][Proof:]{\begin{list}{} 			{\setlength{\leftmargin}{1em} 			\setlength{\rightmargin}{0em}}                         \item {\bf \hspace{-1em}  #1 \ \ }} 
\newcommand{\Expct}[2][]{{\mathbb{E}}_{#1}\left[#2\right]}
\begin{document}

\title{Linear cellular automata, asymptotic randomization, and entropy}

\author{Marcus Pivato\\ {\em Department of Mathematics}, \ {\em Trent University} \\ email: {\tt pivato@xaravve.trentu.ca}}

\maketitle

\begin{abstract} If ${\mathcal{ A}}={{\mathbb{Z}}_{/2}}$, then ${\mathcal{ A}}^{\mathbb{Z}}$ is a compact
abelian group.  A {\bf linear cellular automaton} is a
shift-commuting endomorphism $\Phi$ of ${\mathcal{ A}}^{\mathbb{Z}}$.  If $\mu$ is a
probability measure on ${\mathcal{ A}}^{\mathbb{Z}}$, then $\Phi$ {\bf asymptotically
randomizes} $\mu$ if $\Phi^j \mu$ converges to the Haar measure as
$j{\rightarrow}{\infty}$, for $j$ in a subset of Ces\`aro  density one.  Via
counterexamples, we show that nonzero entropy of $\mu$ is neither
necessary nor sufficient for asymptotic randomization.
\end{abstract}
  If ${\mathcal{ A}}={{\mathbb{Z}}_{/2}}$ (with discrete topology), then
${\mathcal{ A}}^{\mathbb{Z}}$ (with Tychonoff topology) is a compact abelian group.  Let
$ {{{\boldsymbol{\sigma}}}^{}} :{\mathcal{ A}}^{\mathbb{Z}}{{\longrightarrow}}{\mathcal{ A}}^{\mathbb{Z}}$ be the shift map (ie. $ {{{\boldsymbol{\sigma}}}^{}} ({\mathbf{ a}}) \ = \
{\left[a'_z  |_{z\in{\mathbb{Z}}}^{} \right]}$, where $a'_z = a_{z-1}$, \ $\forall
z\in{\mathbb{Z}}$).  A {\bf linear cellular automaton} (LCA) is a topological
group endomorphism $\Phi:{\mathcal{ A}}^{\mathbb{Z}}{{\longrightarrow}}{\mathcal{ A}}^{\mathbb{Z}}$ that commutes with
$ {{{\boldsymbol{\sigma}}}^{}} $.  Let ${\mathcal{ M}}({\mathcal{ A}}^{\mathbb{Z}})$ be the set of Borel probability
measures on ${\mathcal{ A}}^{\mathbb{Z}}$, and let $\eta\in{\mathcal{ M}}({\mathcal{ A}}^{\mathbb{Z}})$ be the Haar
measure.  If $\mu\in{\mathcal{ M}}({\mathcal{ A}}^{\mathbb{Z}})$, we say that $\Phi$ {\bf
asymptotically randomizes} $\mu$ if there is subset ${\mathbb{J}}\subset{\mathbb{N}}$
of Ces\`aro  density one so that
\ $\displaystyle{\mathbf{ w}}\!{\mathbf{ k}}^*\!\!-\!\!\!\lim_{{\mathbb{J}}\ni j{\rightarrow}{\infty}} \Phi^j\mu \ = \  \eta$.

  LCA randomize a broad class of probability measures, including
Bernoulli measures, Markov chains, and Markov random fields
\cite{MaassMartinez,MaassMartinezII,MaassHostMartinez,Lind,
PivatoYassawi2,PivatoYassawi1,FerMaassMartNey}.  One of
the common factors in all these cases is positive entropy.
Conversely,  randomization is impossible for many
zero-entropy measures, such as quasiperiodic or rank one systems
\cite{PivatoYassawi3}.

   Is positive entropy a necessary/sufficient condition for asymptotic
randomization?  We will refute both possibilities, by constructing, in
\S\ref{S:zero.random}, a zero-entropy measure which  asymptotically randomizes,
and in  \S\ref{S:nonzero.norandom}, a positive entropy, ergodic measure which
doesn't.
\paragraph*{Preliminaries:}  If $\Phi$ is any LCA on ${\mathcal{ A}}^{\mathbb{Z}}$,
then there is a finite set ${\mathbb{V}}\subset{\mathbb{Z}}$ 
 so that $\Phi$ can be written as the {\bf polynomial of shifts}
 $\Phi = \sum_{v\in{\mathbb{V}}}  {{{\boldsymbol{\sigma}}}^{v}} $.  This means, for any
${\mathbf{ a}}\in{\mathcal{ A}}^{\mathbb{Z}}$, that $\Phi({\mathbf{ a}})_z =  \sum_{v\in{\mathbb{V}}} a_{z+v}$
for all $z\in{\mathbb{Z}}$.  The advantage of this notation is that composition
of LCA corresponds to multiplication of their respective polynomials.
We can thus apply methods from polynomial algebra over finite fields.

  In particular, we can compute binomial coefficients, mod $2$, as follows:
If $n\in{\mathbb{N}}$, let $\{n^{(i)}\}_{i=0}^{\infty}$ be the
the {\bf binary representation} of $n$, so that
$n \ = \ \sum_{i=0}^{\infty} n^{(i)} 2^i$.  If $N\in{\mathbb{N}}$, then
write ``$n\ll N$'' if $n^{(i)} \leq N^{(i)}$ for all $i\in{\left[ 0..{\infty} \right)}$.
Then {\bf Lucas Theorem} \cite{Lucas} states:
\[
  \left( N\atop n\right) \quad=\quad {\left\{ \begin{array}{rcl}                                  1 &&\mbox{\ if \ } n \ll N \\
					   0 &&\mbox{\ if \ } n \not\ll N                                \end{array}  \right.  }
\quad\pmod{2}
\]
  The simplest nontrivial LCA is the {\bf Ledrappier automaton}
$\Phi=1+ {{{\boldsymbol{\sigma}}}^{}} $. Let ${\mathbb{L}}\left(N\right) = {\left\{ \ell\in{\mathbb{N}} \; ; \; \ell \ll N \right\} }$.
A consequence of Lucas theorem is that
\[
  \Phi^n \quad=\quad \sum_{\ell\in{\mathbb{L}}\left(N\right)}  {{{\boldsymbol{\sigma}}}^{\ell}} .
\]
Thus, the `geometry' of ${\mathbb{L}}\left(N\right)$, as a subset of ${\mathbb{N}}$, determines
the dynamics of $\Phi^N$.

\section{A Zero-Entropy measure that Randomizes
\label{S:zero.random}}

\begin{figure}
{\scriptsize
\[
\begin{array}{cccccccccccccccccccccccccccccccccccccccccccccccccccccc}
 & & & & & & & &  & & & & & & & & p^4_1 &p^4_2 &p^4_3 &p^4_4 \ldots\\\\
 & & & & & & & & p^3_1 &p^3_2 &p^3_3 &p^3_4 &p^3_5 &p^3_6 &p^3_7 &p^3_8 &
 & & &  \\ \\
 & & & &p^2_1&p^2_2&p^2_3&p^2_4& & & & &p^2_1&p^2_2&p^2_3&p^2_4& & & &\\\\ 
 & & p^1_1 & p^1_2 & & & p^1_1 & p^1_2 & & & p^1_1 & p^1_2 & & & p^1_1 & p^1_2 & & & p^1_1 & p^1_2 \ldots \\\\
&p^0_1& &p^0_1& &p^0_1& &p^0_1& &p^0_1& &p^0_1& &p^0_1& &p^0_1& &p^0_1& &p^0_1
\ldots\\
\hline
\ldots a^{\infty}_{1} & a^{\infty}_{2} & a^{\infty}_{3} & a^{\infty}_{4} & a^{\infty}_{5} & a^{\infty}_{6} & a^{\infty}_{7} & a^{\infty}_{8} & a^{\infty}_{9} & a^{\infty}_{10} & a^{\infty}_{11} & a^{\infty}_{12} & a^{\infty}_{13} & a^{\infty}_{14} & a^{\infty}_{15} & a^{\infty}_{16} & a^{\infty}_{17} & a^{\infty}_{18} & a^{\infty}_{19} & a^{\infty}_{20}  \ldots\\
&&&&&\ldots& a_{1} & a_{2} & a_{3} & a_{4} & a_{5} & a_{6} & a_{7} & a_{8} & \ldots
\end{array}
\]}
\caption{The construction of random sequence ${\mathbf{ a}}^{{\infty}}$, 
and the approximation of ${\mathbf{ a}}$ as a random translate of ${\mathbf{ a}}^{{\infty}}$. \label{fig:construct}}
\end{figure} 

  Let $\alpha=\frac{1}{2^{1/5}}$.  For any $n\in{\mathbb{N}}$, let
$\rho_n$ be the probability distribution on ${\mathcal{ A}}={{\mathbb{Z}}_{/2}}$ so that
 \begin{equation} 
\label{prob.def}
  \rho_n\{1\} = \alpha^{n} \hspace{2em} \mbox{and} \ \
 \rho_n\{0\} = 1-\alpha^{n}
 \end{equation} 
  For each $n\in{\mathbb{N}}$, we will construct a random sequence ${\mathbf{ a}}^{n}
\in{\mathcal{ A}}^{\mathbb{Z}}$ as follows.  First, define ${\mathbf{ a}}^{0} =
[\ldots0000\ldots]$.  Now, suppose, inductively, that we have
${\mathbf{ a}}^{n}$.  Let $p^n_1,p^n_2,\ldots,p^n_{2^n}$ be independent
${\mathcal{ A}}$-valued random variables with distribution $\rho_n$.  Let
${\mathbf{ p}}^{n}\in{\mathcal{ A}}^{\mathbb{Z}}$ be the random, $2^{n+1}$-periodic sequence
\[
  {\mathbf{ p}}^{n} \ = \ [\ldots,\underbrace{0,0,\ldots,0}_{2^n},
p^n_1,p^n_2,\ldots,p^n_{2^n},\underbrace{0,0,\ldots,0}_{2^n},p^n_1,p^n_2,\ldots,p^n_{2^n},\ldots],
\]
 and define ${\mathbf{ a}}^{n+1} = {\mathbf{ a}}^{n} + {\mathbf{ p}}^{n}$.

  Let $\mu_n\in{\mathcal{ M}}({\mathcal{ A}}^{\mathbb{Z}})$ be the distribution of ${\mathbf{ a}}^{n}$, and
let $\widetilde\mu_n = \displaystyle \frac{1}{2^n}\sum_{i=1}^{2^n}  {{{\boldsymbol{\sigma}}}^{i}} (\mu_n)$ be
the stationary average of $\mu_n$. Finally, let $\mu = \displaystyle
{\mathbf{ w}}\!{\mathbf{ k}}^*\!\!-\!\!\!\lim_{n{\rightarrow}{\infty}} \widetilde\mu_n$.

  Let $\mu_{{\infty}}$ be the probability distribution of the
random sequence ${\mathbf{ a}}^{{\infty}} = \displaystyle \sum_{n=1}^{\infty} {\mathbf{ p}}^{n}$ (see Figure
\ref{fig:construct}).  Then, $\mu_{\infty} = \displaystyle {\mathbf{ w}}\!{\mathbf{ k}}^*\!\!-\!\!\!\lim_{n{\rightarrow}{\infty}} \mu_n$, and
loosely speaking, $\mu$ is the `$ {{{\boldsymbol{\sigma}}}^{}} $-ergodic average' of $\mu_{{\infty}}$.
Hence, if ${\mathbf{ a}}$ is a $\mu$-random sequence, we can think of ${\mathbf{ a}}$ as
obtained by shifting ${\mathbf{ a}}^{{\infty}}$ by some random amount.

  One way to think of ${\mathbf{ a}}^{{\infty}}$ is as a `randomly generated T\"oplitz
sequence'.  Another way is to imagine ${\mathbf{ a}}^{{\infty}}$ as generated by a process
of `duplication with error'.  Let ${\mathbf{ w}}^{0} = [0]$ be a word of
length 1.  Suppose, inductively, that we have ${\mathbf{ w}}^{n} = [w_1
w_2\ldots w_{2^n}]$.  Let $\widetilde{\mathbf{ w}}^{n} = [{\widetilde{w}}_1
{\widetilde{w}}_2\ldots {\widetilde{w}}_n]$ be an `imperfect copy' of ${\mathbf{ w}}^{n}$: \ \ 
for each $m\in{\left[ 1..2^n \right]}$, \ ${\widetilde{w}}_m = w_m + p^n_m$, where
$p^n_1,p^n_2\ldots,p^n_{2^n}$ are the independent $\rho_n$-distributed
variables from before, which act as `copying errors'.  Let ${\mathbf{ w}}^{n+1}
\ = \ {\mathbf{ w}}^n \widetilde{\mathbf{ w}}^n$.  Then
${\mathbf{ a}}^{{\infty}}$ is the limit of ${\mathbf{ w}}^{n}$ as $n{\rightarrow}{\infty}$.

  If $\nu$ is a probability measure on $\{0,1\}$, then let
\[
  H(\nu) \quad  = \quad  -\nu\{0\}\log(\nu\{0\})  \ - \ \nu\{1\}\log(\nu\{1\})
\]
  be the entropy of $\nu$.  If $b$ is a $\nu$-random variable, and
  $\chi:{\mathcal{ A}}{{\longrightarrow}}{\mathbb{R}}$ is some function, then let
\[
  \Expct[\nu]{\chi(b)} \quad = \quad  \nu\{0\}\chi(0) \ + \  \nu\{1\}\chi(1)
\]
  be the expected value of $\chi(b)$.

\begin{lemma}\label{perturb.lemma}  
 Let ${\mathbf{ a}}\in{\mathcal{ A}}^{\mathbb{Z}}$ be a $\mu$-random sequence, and fix $n\in{\mathbb{N}}$.
 Then for all $m\in{\left[ 1..2^n \right]}$,
 $a_{m+2^n} = a_{m} + d_m$, where $d_1,\ldots,d_{2^n}$ are independent
random variables with distributions $\delta_1,\ldots,\delta_{2^n}$, such that,
for all $m$,

 {\rm{\bf (i)}} \quad \ If $n>5$, then \ 
 $\alpha^n \ < \ \delta_m\{1\} \ < \ 8\alpha^{n}$;

 {\rm{\bf (ii)}} \quad If $n>20$, then \  
$H(\delta_m) \ \ < \ \  2 n\cdot\alpha^{n}$.
 \end{lemma}
\bprf  By construction, there is some $k\in{\mathbb{Z}}$ so that
${\mathbf{ a}}$ looks like $ {{{\boldsymbol{\sigma}}}^{k}} ({\mathbf{ a}}^{{\infty}})$ in a neighbourhood
around $0$.  To be precise,
$a_{m} = a^{{\infty}}_{k+m}$ for all $m\in{\left[ 1..2^{n+1} \right]}$.  

  For example, consider Figure \ref{fig:construct}, and let $n=2$,
so that $2^n=4$; suppose $k=6$.  Thus,
\[
 [a_1,a_2,\ldots,a_8] 
\ = \ 
 [a^{\infty}_7,a^{\infty}_8,\ldots,a^{\infty}_{14}]
\]
  Thus, $d_1 = a_5-a_1 \ = \ a^{\infty}_{11} - a^{\infty}_7 \ = \ 
p^3_3 - p^2_3  =  p^3_3 + p^2_3$.  Similarly, for any $m\in{\left[ 1..2^n \right]}$,
\[
 d_m \ \  =  \ \  a_{m+2^n} - a_m \ \ = \  \  a^{\infty}_{k+m+2^n} - a^{\infty}_{k+m} 
 \ \ = \ \  p^n_{m_0} +  p^{n_1}_{m_1} + \ldots +  p^{n_J}_{m_J}
\]
  where $n = n_0< n_1 < \ldots < n_J$ (and depend on $m$ and $k$), and where
$p^{n_j}_{m_j}$ has distribution $\rho_{n_j}$.  For
all $j\in{\left[ 0..J \right]}$, let 
$P_j={\sf Prob}\left( \ \rule[-0.5em]{0em}{1em}       \begin{minipage}{40em}       \begin{tabbing}         $\displaystyle\sum_{i=j}^J p^{n_i}_{m_i}$ \ is odd        \end{tabbing}      \end{minipage} \ \right)$.
Thus, 
\begin{eqnarray*}
 \delta_m\{1\} 
&=&
P_0
\quad=\quad
\rho_{n}\{0\}\cdot P_1 \ + \ \rho_{n}\{1\}\cdot (1-P_1) 
\quad=\quad 
(1-\alpha^n) \cdot P_1 \ + \ \alpha^n \cdot (1-P_1)\\
&=&
\alpha^n \ + \ (1-2\alpha^n)\cdot P_1
\quad\geq\quad
\alpha^n, \quad
\mbox{(because $n>5$, so  $\alpha^n<\frac{1}{2}$, so $1-2\alpha^n>0$.)}
\end{eqnarray*}
  This holds for any $k\in{\mathbb{Z}}$.  Average over all $k$ to get the
lower bound in {\bf(i)}.

Also, for any $j\in{\left[ 1..J \right)}$, we have
\[
P_j
\quad=\quad
(1-\alpha^{n_j}) \cdot P_{j+1} \ + \ \alpha^{n_j} \cdot (1-P_{j+1})\\
\quad=\quad 
P_{j+1} \ + \ (1-2P_{j+1})\alpha^{n_j}
\quad\leq\quad
P_{j+1} \ + \ \alpha^{n_j},
\]
and $P_J=\alpha^{n_J}$.  Hence, inductively,
\[
 \delta_m\{1\} \quad=\quad  P_0
 \quad\leq\quad
 \alpha^{n_0} + \alpha^{n_1} + \ldots + \alpha^{n_J}
 \quad\leq\quad
  \sum_{i=n_0}^{\infty} \alpha^i
 \quad = \quad
 \alpha^{n_0} \frac{1}{1-\alpha}
 \quad < \quad
 8\alpha^{n_0}.
\]
Again, this holds for any
$k\in{\mathbb{Z}}$; average over all $k$ to get the
upper bound in {\bf(i)}.

\paragraph{Proof of {\bf(ii)}:}
\quad
  If $\nu\{1\}<\frac{1}{2}$, then $H(\nu)$ decreases as $\nu\{1\}$ decreases.
  If $n>20$, then Part {\bf(i)} says
$\delta_m\{1\} \ \leq \ 8\alpha^n  \ = \ 2^{3-n/5} \ \leq \ \frac{1}{2}$; \  
hence,
\begin{eqnarray*}
H(\delta_m) & \leq & 
  -8\alpha^{n}\log_2 \left(8\alpha^{n}\right) \ - \
 \left(1-8\alpha^{n}\right)\log_2 \left(1-8\alpha^{n}\right)
\\& <&
  8 \alpha^{n}\left(\frac{n}{5}-3\right) \ + \ \left(1-8\alpha^{n}\right)\cdot 
\overbrace{2\cdot 8\alpha^{n}}^{(a)}
 \quad = \quad
 8  \left(\frac{n}{5}-1-16\alpha^n\right)\cdot \alpha^{n}\\
 & < & \frac{8}{5}n\cdot \alpha^{n} 
\quad<\quad 2n\cdot \alpha^n,
\end{eqnarray*}
where $(a)$ is because, for small $\epsilon$, $\log(1-\epsilon) \approx -\epsilon$,
thus,  $-\log(1-\epsilon) < 2\epsilon$.
 {\tt \hrulefill $\Box$ } \end{list}  \medskip  

\begin{lemma}\label{entropy.lemma}   
  $h(\mu) = 0$.
 \end{lemma}
\bprf
  Suppose ${\mathbf{ a}}\in{\mathcal{ A}}^{\mathbb{Z}}$ is a $\mu$-random sequence.  Fix $n>20$;
we want to compute the conditional entropy 
$H\left({\mathbf{ a}}\raisebox{-0.3em}{$\left|_{{\left( 2^n..2^{n+1} \right]}}\right.$} | {\mathbf{ a}}\raisebox{-0.3em}{$\left|_{{\left[ 1..2^n \right]}}\right.$}\right)$.  By Lemma
\ref{perturb.lemma}{\bf(ii)}, we know that, for all $m\in{\left[ 1..2^n \right]}$,
\quad $a_{2^n+m} \ = \ a_m + d_m$, \ 
 where $d_1,d_2,\ldots,d_{2^n}$ are independent random variables
with distributions $\delta_1,\ldots,\delta_{2^n}$,
such that $H(\delta_m) \ < \ 2n\alpha^n$.  Thus,
\[
H\left({\mathbf{ a}}\raisebox{-0.3em}{$\left|_{{\left( 2^n..2^{n+1} \right]}}\right.$} \right|\left. {\mathbf{ a}}\raisebox{-0.3em}{$\left|_{{\left[ 1..2^n \right]}}\right.$}\right)
\ =\quad 
H\left(d_1,d_2,\ldots,d_{2^n}\right) 
\ = \quad
\sum_{m=1}^{2^n} H(\delta_m)
\ < \quad
 2^n\cdot 2 n \alpha^n
\ = \ 2 n \cdot (2\alpha)^n.
\]
Thus, for any $N>20$,
\begin{eqnarray*}
\lefteqn{
H\left({\mathbf{ a}}\raisebox{-0.3em}{$\left|_{{\left[ 1..2^N \right]}}\right.$} \right|\left. {\mathbf{ a}}\raisebox{-0.3em}{$\left|_{{\left[ 1..2^{20} \right]}}\right.$} \right)
\quad = \quad
\sum_{n=20}^{N-1} 
H\left({\mathbf{ a}}\raisebox{-0.3em}{$\left|_{{\left( 2^n..2^{n+1} \right]}}\right.$}  \right|\left. {\mathbf{ a}}\raisebox{-0.3em}{$\left|_{{\left[ 1..2^n \right]}}\right.$}\right) 
\quad < \quad
\sum_{n=20}^{N-1}
   2 n \cdot (2\alpha)^n} \\
& \hspace{4em}\leq &
 2 N \cdot (2\alpha)^{20}\sum_{n=0}^{N-21}
  (2\alpha)^n 
\quad = \quad
 2 N  \cdot (2\alpha)^{20} \frac{(2\alpha)^{N-20}-1}{2\alpha -1}
\quad \leq \quad
 c N \cdot (2\alpha)^N,
\end{eqnarray*}
where $c$ is a constant.  Thus, if $H_0 = 
\left({\mathbf{ a}}\raisebox{-0.3em}{$\left|_{{\left[ 1..2^{20} \right]}}\right.$} \right)$, then
\[ 
H\left({\mathbf{ a}}\raisebox{-0.3em}{$\left|_{{\left[ 1..2^N \right]}}\right.$}  \right)
\quad=\quad
H\left({\mathbf{ a}}\raisebox{-0.3em}{$\left|_{{\left[ 1..2^N \right]}}\right.$} \right|\left. {\mathbf{ a}}\raisebox{-0.3em}{$\left|_{{\left[ 1..2^{20} \right]}}\right.$} \right) \ + \  H_0
\quad \leq \quad
 c N \cdot (2\alpha)^N \ + \ H_0.
\]
\begin{eqnarray*}
\mbox{Hence,}\quad\quad h(\mu)
&=&
\lim_{M{\rightarrow}{\infty}}\ \frac{1}{M}\ H\left({\mathbf{ a}}\raisebox{-0.3em}{$\left|_{{\left[ 1..M \right]}}\right.$} \right)
\quad= \quad
\lim_{N{\rightarrow}{\infty}}\ \frac{1}{2^N}\ H\left({\mathbf{ a}}\raisebox{-0.3em}{$\left|_{{\left[ 1..2^N \right]}}\right.$} \right)\\
&=&
\lim_{N{\rightarrow}{\infty}}\ \frac{c N \cdot (2\alpha)^N + H_0 }{2^N} 
\quad\leq \quad 
  c \lim_{N{\rightarrow}{\infty}}\  N \alpha^{N} 
\quad= \quad 0,
\end{eqnarray*}
because $\alpha<1$.
 {\tt \hrulefill $\Box$ } \end{list}  \medskip  

 ${\mathcal{ A}}^{\mathbb{Z}}$ is a compact abelian group; \ let $\widehat{{\mathcal{ A}}^{\mathbb{Z}}}$
 be its group of characters.  The only nontrivial character of ${\mathcal{ A}}={{\mathbb{Z}}_{/2}}$ is the map ${\mathcal{ E}}:{\mathcal{ A}}{{\longrightarrow}}\{\pm1\}$ defined: ${\mathcal{ E}}(a) = (-1)^a$.
If ${{{\mathsf{ 1\!\!1}}}_{{}}}\neq{\boldsymbol{\chi }}\in\widehat{{\mathcal{ A}}^{\mathbb{Z}}}$, then there is a finite subset
${\mathbb{K}}\subset{\mathbb{Z}}$ so that, for any ${\mathbf{ a}}\in{\mathcal{ A}}^{\mathbb{Z}}$,
\ \  ${\boldsymbol{\chi }}({\mathbf{ a}}) \ = \displaystyle \prod_{k\in{\mathbb{K}}} {\mathcal{ E}}(a_k)$.
For all $k\in{\mathbb{K}}$, define $\chi_k:{\mathcal{ A}}^{\mathbb{Z}}{{\longrightarrow}}\{\pm1\}$ by
$\chi_k({\mathbf{ a}}) = {\mathcal{ E}}(a_k)$.  We define ${{\sf rank}\left[{\boldsymbol{\chi }}\right]} \quad=\quad {\sf card}\left[{\mathbb{K}}\right]$,
\ and  \ 
${\sf diam}\left[{\boldsymbol{\chi }}\right] \ = \ \max({\mathbb{K}}) - \min({\mathbb{K}}) + 1$, and
write
\ ``$\displaystyle {\boldsymbol{\chi }} \ = \ \bigotimes_{k\in{\mathbb{K}}} \chi_k$''.

\begin{lemma}\label{asympt.rand.lemma}   
Let $\Phi=1+ {{{\boldsymbol{\sigma}}}^{}} $.  Then $\Phi$ asymptotically randomizes $\mu$.
 \end{lemma}
\bprf
As in the proof of \cite[Theorem 12]{PivatoYassawi1}, it is sufficient
to prove that, for any ${{{\mathsf{ 1\!\!1}}}_{{}}}\neq{\boldsymbol{\chi }}\in\widehat{{\mathcal{ A}}^{\mathbb{Z}}}$,
there is a set ${\mathbb{J}}\subset{\mathbb{N}}$ of density one so that
$\displaystyle \lim_{{\mathbb{J}}\ni j {\rightarrow}{\infty}} {\left\langle \Phi^j\mu, \ {\boldsymbol{\chi }} \right\rangle } \quad= \quad0$.
Let $\displaystyle {\boldsymbol{\chi }} \ = \ \bigotimes_{k\in{\mathbb{K}}} \chi_k$.  Let $N\geq
\log_2({\sf diam}\left[{\boldsymbol{\chi }}\right])$, and 
assume without loss of generality that ${\mathbb{K}}\subset{\left[ 0..2^N \right]}$.

  Let $n\in{\mathbb{N}}$ be some large integer with binary expansion
$\{n^{(i)}\}_{i=0}^I$, where
$I=\lfloor\log_2(n)\rfloor$.  For generic $n\in{\mathbb{N}}$
(ie. for a set of $n$ of Ces\`aro  density 1), we can find  $J\in{\mathbb{N}}$
 so that
\[
  \mbox{\bf(G1)} \quad N+2 < \ J \ < I/2;
\hspace{5em}\mbox{\ and \ }\hspace{4em}
 \mbox{\bf(G2)} \quad n^{(J-2)} = n^{(J-1)} = 0.
\]
  Thus, $n \ = \ n_0 \ + \ 2^J n_1$, where
$\displaystyle
n_0 \ = \ \sum_{j=0}^{J-3} n^{(j)} 2^j$
\quad and \quad
$\displaystyle n_1 \ = \ \sum_{j=J}^{I} n^{(j)} 2^{j-J}$.

  Thus, ${\mathbb{L}}\left(n\right) \ = \ {\mathbb{L}}\left(n_0\right) \ + \ 2^J\cdot {\mathbb{L}}\left(n_1\right)$.
Now, let \ $\displaystyle
  {\boldsymbol{\xi }}_0 \ = \ \bigotimes_{\ell\in{\mathbb{L}}\left(n_0\right)} {\boldsymbol{\chi }}\circ {{{\boldsymbol{\sigma}}}^{\ell}} $.
Then
\[
  {\boldsymbol{\chi }}\circ\Phi^{n}
\quad= \quad
 \bigotimes_{\ell\in{\mathbb{L}}\left(n\right)} \ {\boldsymbol{\chi }}\circ {{{\boldsymbol{\sigma}}}^{\ell}} 
\quad= \  \quad
 \bigotimes_{\ell_1\in{\mathbb{L}}\left(n_1\right)}  \  \bigotimes_{\ell_0\in{\mathbb{L}}\left(n_0\right)} \
 {\boldsymbol{\chi }}\circ {{{\boldsymbol{\sigma}}}^{\ell_0 + 2^J\ell_1}} 
\quad= \quad
\bigotimes_{\ell_1\in{\mathbb{L}}\left(n_1\right)} \
 {\boldsymbol{\xi }}_0 \circ  {{{\boldsymbol{\sigma}}}^{2^J \ell_1}} 
\]
is a `product of translates'  of ${\boldsymbol{\xi }}_0$.
 These translates do not overlap, because
\[
{\sf diam}\left[{\boldsymbol{\xi }}_0\right] \quad\leq \quad{\sf diam}\left[{\boldsymbol{\chi }}\right] + n_0 \quad\leq \quad2^N + 2^{J-2} 
\quad<_{\mathrm{by} \ \mathbf{(G1)}} \quad2^{J-2} + 2^{J-2} 
\quad< \quad2^J.
\]
  Let ${\mathbb{I}} = {\left\{ j\in{\left[ J..I \right]} \; ; \; n^{(j)} =1 \right\} }$.
Let $\epsilon>0$ be small; \  then for generic $n\in{\mathbb{N}}$, we can assume
\begin{itemize}
\item[{\bf (G3)}] \hspace{2em} ${\sf card}\left[{\mathbb{I}}\right] \  \ \geq \  \  \frac{1}{2} (I-J) - \epsilon$.
\end{itemize}
Since  $\alpha=\frac{1}{2^{1/5}}$, we can
find $\beta$ such that $\frac{1}{\alpha} \ < \ \beta \ < \ 2^{1/4}$. 
Thus, if $M={\sf card}\left[{\mathbb{I}}\right]-1$, then 
 \begin{equation} 
\label{M.bigger.beta}
M \quad\geq_{\mathrm{by} \ \mathbf{(G3)}}
 \quad\frac{1}{2} (I-J) - \epsilon - 1
\quad >_{\mathrm{by} \ \mathbf{(G1)}} \quad \frac{1}{4}I - \epsilon-1 \quad>
 \quad\log_2(\beta) I,
 \end{equation} 
 because $\log_2(\beta) \ < \ \frac{1}{4}$ and $I$ is large,
while $\epsilon$ is small.

  Suppose ${\mathbb{I}} = \{ i_1 < i_2 < \ldots < i_{M+1} = I\}$.  For 
each $m\in{\left[ 0..M \right]}$, define \
 ${\boldsymbol{\xi }}_{m+1} =  {\boldsymbol{\xi }}_m  \otimes \ \left({\boldsymbol{\xi }}_m\circ {{{\boldsymbol{\sigma}}}^{L}} \right)$, \
where $L=2^{i_m}$. \ 
Thus, \  ${\boldsymbol{\chi }}\circ \Phi^n  \ =  \ {\boldsymbol{\xi }}_{M+1}$.  

  Let $r = {{\sf rank}\left[{\boldsymbol{\xi }}_0\right]}$.  Then for all $m\in{\left[ 1..M+1 \right]}$, \ \ 
${{\sf rank}\left[{\boldsymbol{\xi }}_m\right]} \ = \ 2^m \cdot r$.  In particular, define
 \begin{equation} 
 R \quad = \quad {{\sf rank}\left[{\boldsymbol{\xi }}_M\right]} \quad = \quad  2^M \cdot r  
\quad >_{\mathrm{by} \ (\ref{M.bigger.beta})} \quad
 \beta^I \cdot r. 
\label{R.defn}
 \end{equation} 
Thus, \
$\displaystyle  {\boldsymbol{\xi }}_M  \ =  \ \bigotimes_{x\in{\mathbb{X}}} \xi_x$, \ 
where ${\mathbb{X}}\subset{\mathbb{Z}}$ is a subset with ${\sf card}\left[{\mathbb{X}}\right]=R$.  Thus,
if ${\mathbf{ a}}\in{\mathcal{ A}}^{\mathbb{Z}}$ is a $\mu$-random sequence, then
\[
{\boldsymbol{\xi }}_{M+1}({\mathbf{ a}}) \ = \ 
{\boldsymbol{\xi }}_{M}({\mathbf{ a}}) \ \cdot \  \left({\boldsymbol{\xi }}_M\circ {{{\boldsymbol{\sigma}}}^{2^I}} ({\mathbf{ a}})\right)
\ = \  \prod_{x\in{\mathbb{X}}} \xi_x(a_x) \cdot  \xi_x \left(a_{x+2^I}\right)
\ = \  \prod_{x\in{\mathbb{X}}} \xi_x\left(a_x + a_{x+2^I}\right)
\ = \ \prod_{x\in{\mathbb{X}}} \xi_x\left(d_x\right),
\]
where $\{d_x\}_{x\in{\mathbb{X}}}$ are  independent random variables
as in Lemma \ref{perturb.lemma}. 
 Let $d_x$ have distribution $\delta_x$;
then $\delta_x\{1\} \geq \alpha^I$, by Lemma \ref{perturb.lemma}{\bf(i)}.  Thus,
\[
 \Expct[\delta_x]{ \xi_x\left(d_x\right)} 
\quad=\quad
\delta_x\{0\}-\delta_x\{1\} 
\quad=\quad 
1-2\delta_x\{1\}  
 \quad \leq \quad 
1-2\cdot\alpha^{I}
\quad = \quad 
\frac{2\alpha^{-I} - 1}{2\alpha^{-I}}. 
\]
\[
\mbox{Thus,}\quad
{\left\langle \mu, \ \ {\boldsymbol{\chi }}\circ\Phi^n \right\rangle }
\quad= \quad
\Expct{\prod_{x\in{\mathbb{X}}} \xi_x\left(d_x\right)}
\quad= \quad
\prod_{x\in{\mathbb{X}}} \Expct[\delta_x]{ \xi_x\left(d_x\right)}
\quad\leq \quad
\left(\frac{2\alpha^{-I} - 1}{2\alpha^{-I}} \right)^R.
\]
\begin{eqnarray*}
\mbox{Thus,}\quad
\log \left|\rule[-0.5em]{0em}{1em} {\left\langle \mu, \ \ {\boldsymbol{\chi }}\circ\Phi^n \right\rangle }\right|
& \leq & 
 R \cdot \left[\rule[-0.5em]{0em}{1em} 
\log\left(2\alpha^{-I} - 1\right)-\log(2\alpha^{-I})\right] 
\quad \leq_{(\ast)} \quad
 -R \cdot \log'\left(2\alpha^{-I}\right)\\
& = & \frac{-R}{2\alpha^{-I}}
\quad <_{\mathrm{by} \ (\ref{R.defn})} \quad \frac{-\beta^I\, r}{2\alpha^{-I}} 
\quad = \quad
- \frac{r}{2} \ (\alpha\beta)^I,
\end{eqnarray*}
$(\ast)$ is because $\log$ is a decreasing function.
But $\beta> \frac{1}{\alpha}$, so  $\alpha\beta>1$.
Thus, if ${\mathbb{J}}\subset{\mathbb{N}}$ is
the set of all $n\in{\mathbb{N}}$ satisfying the generic hypotheses {\bf(G1-G3)}, then
$\displaystyle \lim_{{\mathbb{J}}\ni n{\rightarrow}{\infty}} \log  \left|\rule[-0.5em]{0em}{1em} {\left\langle \mu, \  {\boldsymbol{\chi }}\circ\Phi^n \right\rangle }\right|
\ = \ - \frac{r}{2} \ \lim_{I{\rightarrow}{\infty}} \ 
(\alpha\beta)^I \ = \ -{\infty}$.  Hence $\displaystyle \lim_{{\mathbb{J}}\ni
n{\rightarrow}{\infty}} \left|{\left\langle \mu, \ {\boldsymbol{\chi }}\circ\Phi^n \right\rangle }\right|
\  = \ 0$.
 {\tt \hrulefill $\Box$ } \end{list}  \medskip

\section{\label{S:nonzero.norandom}
A nonzero-entropy measure that doesn't randomize}

  Let $Q=2^k$ for some $k>0$.  Treat ${\mathcal{ A}}$ as a field, and
${\mathcal{ A}}^{\left[ 1..Q \right]}$ as a $Q$-dimensional vector space over ${\mathcal{ A}}$.  Let
$R<Q$, and let 
${\mathcal{ B}}\subset{\mathcal{ A}}^{\left[ 1..Q \right]}$ be an $R$-dimensional vector subspace.
Define $\psi:{\mathcal{ B}}^{\mathbb{Z}}{{\longrightarrow}}{\mathcal{ A}}^{\mathbb{Z}}$ as follows: if
${\mathbf{ B}}\in{\mathcal{ B}}^{\mathbb{Z}}$, where ${\mathbf{ B}} = [{\mathbf{ b}}^{(z)}]_{z\in{\mathbb{Z}}}$ and
${\mathbf{ b}}^{(z)} \ = \ [b^{(z)}_q]_{q=1}^{Q}$, then
\[
  \psi({\mathbf{ B}}) \quad=\quad
\left[ \ldots b^{(-1)}_1,b^{(-1)}_2,\ldots,b^{(-1)}_Q, \
b^{(0)}_1,b^{(0)}_2,\ldots,b^{(0)}_Q, \
b^{(1)}_1,b^{(1)}_2,\ldots,b^{(1)}_Q, \ldots\right]
\]
where $\psi({\mathbf{ B}})_1 \ = \ b^{(0)}_1$.  Let ${\mathfrak{ B}} = {\sf image}\left[\psi\right]\subset
{\mathcal{ A}}^{\mathbb{Z}}$, with orbit closure ${\mathfrak{ X}} \ = \  \bigsqcup_{q=1}^Q  {{{\boldsymbol{\sigma}}}^{q}} ({\mathfrak{ B}})$.

\begin{lemma}\label{L:nonzero.nonrandom.0}  
Let $\Phi=1+ {{{\boldsymbol{\sigma}}}^{}} $.  Then for any $n\in{\mathbb{N}}$, \quad
 $\Phi^{nQ} ({\mathfrak{ X}})\subset{\mathfrak{ X}}$.
 \end{lemma}
\bprf  Lucas' theorem implies that $\Phi^{Q} = 1+ {{{\boldsymbol{\sigma}}}^{Q}} $;
We claim that $\Phi^{Q}({\mathfrak{ B}})\subset{\mathfrak{ B}}$.  To see this,
let ${\mathbf{ a}}=\psi({\mathbf{ B}})$ for some ${\mathbf{ B}} = [{\mathbf{ b}}^{(z)}]_{z\in{\mathbb{Z}}}
\in{\mathcal{ B}}^{\mathbb{Z}}$. Then
$\Phi^{Q}({\mathbf{ a}}) \ = \  {\mathbf{ a}}+ {{{\boldsymbol{\sigma}}}^{Q}} ({\mathbf{ a}})
\ = \ 
\psi\left(\rule[-0.5em]{0em}{1em}{\mathbf{ b}} +  {{{\boldsymbol{\sigma}}}^{}} ({\mathbf{ b}})\right)$,
so $\Phi^{Q}({\mathbf{ a}})\in{\mathfrak{ B}}$ also.  
  Hence, $\Phi^{Q} ({\mathfrak{ X}})\subset{\mathfrak{ X}}$; hence $\Phi^{nQ} ({\mathfrak{ X}})\subset{\mathfrak{ X}}$
for all $n\in{\mathbb{N}}$.
 {\tt \hrulefill $\Box$ } \end{list}  \medskip  

  Let $\nu\in{\mathcal{ M}}({\mathcal{ B}}^{\mathbb{Z}})$ be the uniformly distributed Bernoulli
measure; let ${\widetilde{\nu}}=\phi(\nu)$, and let $\mu \ = \  \sum_{q=1}^Q
 {{{\boldsymbol{\sigma}}}^{q}} ({\widetilde{\nu}})$.  Then $\mu$ is a $ {{{\boldsymbol{\sigma}}}^{}} $-ergodic
measure, and ${\sf supp}\left(\mu\right)={\mathfrak{ X}}$.

\begin{lemma}\label{L:nonzero.nonrandom.1}  $h(\mu) \ =\  \frac{R}{Q}\log_2(P)$. \end{lemma}
\bprf Every $Q$ symbols of an element of ${\mathfrak{ X}}$ corresponds to a
 single symbol of some element of  ${\mathcal{ B}}^{\mathbb{Z}}$, and $h(\nu)=R\cdot\log_2(P)$.
 {\tt \hrulefill $\Box$ } \end{list}  \medskip  
\begin{lemma}\label{L:nonzero.nonrandom.2}  
 $\Phi$ cannot asymptotically randomize $\mu$.
 \end{lemma}
\bprf  Let $\mu_n=\Phi^n(\mu)$ for all $n\in{\mathbb{N}}$.
Lemma \ref{L:nonzero.nonrandom.0} implies that ${\sf supp}\left(\mu_{nQ}\right)
\subset {\mathfrak{ X}}$ for any $n\in{\mathbb{N}}$.  Thus, the sequence
$\{\mu_n\}_{n=1}^{\infty}$ cannot converge to $\eta$ along a set of
density one.
 {\tt \hrulefill $\Box$ } \end{list}  \medskip

{\footnotesize
\bibliographystyle{plain}
\bibliography{bibliography}

\begin{thebibliography}{1}

\bibitem{MaassMartinez}
{Alejandro Maass} and {Servet Mart{\'i}nez}.
\newblock On {Ces{\`a}ro} limit distribution of a class of permutative cellular
  automata.
\newblock {\em Journal of Statistical Physics}, 90(1-2):435--452, 1998.

\bibitem{MaassMartinezII}
{Alejandro Maass} and {Servet Mart{\'i}nez}.
\newblock Time averages for some classes of expansive one-dimensional cellular
  automata.
\newblock In {Eric Goles} and {Servet Martinez}, editors, {\em Cellular
  Automata \& Complex Systems}, pages 37--54. Kluwer Academic Publishers,
  Dordrecht, 1999.

\bibitem{MaassHostMartinez}
{Servet Martinez} {Alejandro Maass}, {Bernard Host}.
\newblock Uniform {B}ernoulli measure in dynamics of permutative cellular
  automata with algebraic local rules.
\newblock {\em {\rm submitted to} Discrete \& Continuous Dyn. Sys.}, 2002.

\bibitem{Lind}
Doug Lind.
\newblock Applications of ergodic theory and sofic systems to cellular
  automata.
\newblock {\em Physica D}, 10:36--44, 1984.

\bibitem{Lucas}
E.~Lucas.
\newblock Sur les congruences des nombres {Eul\'eriens} et des coefficients
  diff{\'e}rentiels des fonctions trigonom{\'e}triques, suivant un module
  premier.
\newblock {\em Bulletin de la Soc. Math. de France}, 6:49--54, 1878.

\bibitem{PivatoYassawi2}
{Marcus Pivato} and {Reem Yassawi}.
\newblock Limit measures for affine cellular automata {II}.
\newblock {\em {\rm Submitted to} Ergodic Theory \& Dynamical Systems; {\rm
  preprint available:} {{\tt http://arXiv.org/abs/math.DS/0108083}}}, 2001.

\bibitem{PivatoYassawi3}
{Marcus Pivato} and {Reem Yassawi}.
\newblock Asymptotic behaviour of measures with long range correlations under
  the action of cellular automata.
\newblock {\em (in press)}, 2002.
\newblock {\rm preprint available:} {\tt http://arXiv.org/abs/math.DS/0210232}.

\bibitem{PivatoYassawi1}
{Marcus Pivato} and {Reem Yassawi}.
\newblock Limit measures for affine cellular automata.
\newblock {\em Ergodic Theory \& Dynamical Systems}, 22(4):1269--1287, August
  2002.
\newblock ({\tt http://arXiv.org/abs/math.DS/0108082}).

\bibitem{FerMaassMartNey}
{S. Mart{\'i}nez} {P. Ferrari}, {A. Maass} and {P. Ney}.
\newblock Ces{\`a}ro mean distribution of group automata starting from measures
  with summable decay.
\newblock {\em Ergodic Theory \& Dynamical Systems}, 20(6):1657--1670, 2000.

\end{thebibliography}
}

\end{document}